\documentclass[11pt]{article}

\title{K-groups of the quantum homogeneous space $SU_{q}(n)/SU_{q}(n-2)$}
\author{Partha Sarathi Chakraborty and S.Sundar}
\usepackage{latexsym}
\usepackage{amsfonts}
\usepackage{amssymb}
\usepackage{amsmath}
\usepackage{mathrsfs}
\input xypic
\usepackage{hyperref}
\usepackage{tikz}
\usetikzlibrary{matrix,arrows}

\newcommand{\be}{\begin{equation}}
\newcommand{\ee}{\end{equation}}
\newcommand{\bea}{\begin{eqnarray}}
\newcommand{\eea}{\end{eqnarray}}
\newcommand{\bean}{\begin{eqnarray*}}
\newcommand{\eean}{\end{eqnarray*}}
\newcommand{\brray}{\begin{array}}
\newcommand{\erray}{\end{array}}

\newtheorem{dfn}{Definition}[section]
\newtheorem{thm}[dfn]{Theorem}
\newtheorem{lmma}[dfn]{Lemma}
\newtheorem{ppsn}[dfn]{Proposition}
\newtheorem{crl}[dfn]{Corollary}
\newtheorem{xmpl}[dfn]{Example}
\newtheorem{rmrk}[dfn]{Remark}

\newcommand{\bdfn}{\begin{dfn}\rm}
\newcommand{\bthm}{\begin{thm}}
\newcommand{\blmma}{\begin{lmma}}
\newcommand{\bppsn}{\begin{ppsn}}
\newcommand{\bcrlre}{\begin{crlre}}
\newcommand{\bxmpl}{\begin{xmpl}}
\newcommand{\brmrk}{\begin{rmrk}\rm}

\newcommand{\edfn}{\end{dfn}}
\newcommand{\ethm}{\end{thm}}
\newcommand{\elmma}{\end{lmma}}
\newcommand{\eppsn}{\end{ppsn}}
\newcommand{\ecrlre}{\end{crlre}}
\newcommand{\exmpl}{\end{xmpl}}
\newcommand{\ermrk}{\end{rmrk}}

\begin{document}
\maketitle
\begin{abstract}
Quantum Steiffel manifolds were introduced by Vainerman and Podkolzin in \cite{VP}. They classified the irreducible representations of their underlying $C^*$-algebras. Here we compute   the K groups of the quantum homogeneous spaces
$SU_{q}(n)/SU_{q}(n-2), n\ge 3$.  Specializing to the case $n=3$ we show that the fundamental unitary for quantum $SU(3)$ is nontrivial and forms part of a generating set in the $K_1$. 
\end{abstract}
\section{Introduction}
Quantization of mathematical theories is a major theme of research today. The theory of Quantum Groups and Noncommutative Geometry are two prime examples in this program. Both these programs started in the early eighties. In the setting of operator algebras the theory of quantum groups was initiated independently  by Woronowicz (\cite{Wor_twisted}); Vaksman and Soibelman (\cite{Vak-Soi}). They studied the case of quantum $SU(2)$. Later Woronowicz studied the family of compact quantum groups and obtained Tannaka type duality theorems (\cite{Wor_Tannaka}). Soon the notion of quantum subgroups and quantum homogeneous spaces followed (\cite{Pod}). 

The  Noncommutative Differential Geometry program of Alain Connes also started in the eighties (\cite{ConIHES}). In his interpretation geometric data is encoded in elliptic operators or more generally in specific unbounded K-cycles, which he called spectral triples. It is natural to expect that for compact quantum groups and their homogeneous spaces there should be associated canonical spectral triples. Chakraborty and Pal showed that (\cite{PsPal_SU_q_2})  indeed that is the case for quantum $SU(2)$. In fact  for odd dimensional quantum spheres one can construct finitely summable spectral triples that witnesses Poincare duality (\cite{CP_Tran}). A natural question in this connection is, are these examples somewhat singular or in general one can construct finitely summable spectral triples with further properties like Poincare duality on quantum groups associated with Lie groups or their homogeneous spaces. Even though there are suggestions to construct such spectral triples (\cite{NT_arxiv}) one does not know their nontriviality as a K-cycle. In fact there are suggestions that for quantum groups and their homogeneous spaces one should look for a type III formulation  of noncommutative geometry.  On this formulation also there are currently two
points of view, that of Alain Connes and Henri Moscovici
(\cite{CM_typeIII}), and that of
Carey-Phillips-Rennie (\cite{CPR_arxiv}). Therefore to understand the true nature of interplay between noncommutative geometry and quantum homogeneous spaces it makes sense to take a closer look at these algebras. The underlying $C^*$-algebras of these compact quantum groups were analysed by Soibelman (\cite{So1}) (also \cite{Le-So}) who described the irreducible representations of these algebras. Exploiting their findings Sheu went on to obtain composition sequences for these algebras. He initially obtained the results for $SU_q(3)$ (\cite{Sh-pacific}) and later extended them for the general $SU_q(n)$ (\cite{Sh1}). In this hierarchy of exploration the next thing to look for, would be the K-groups and that is what we are looking for.  But instead of concentrating on the quantum groups we consider quantum analogs of the Steiffel manifolds $SU(n)/SU(n-m)$, introduced by Podkoljin and Vainerman (\cite{VP}). They have already described the structure of irreducible representations of quantum Steiffel manifolds $SU_q(n)/SU_q(n-m)$. We take up the case of $SU_q(n)/SU_q(n-2), n \ge 3$. We obtain the composition sequences for these algebras and then utilising them we compute the K-groups. More importantly as we remarked earlier applications towards noncommutative geometry requires explicit understanding of generators for these K-groups and during our calculation we also achieve that. Specializing to the case of $n=3$ we get the K-groups of quantum $SU(3)$. Here we should remark that probably these K-groups can also be computed using the variant of KK-theory introduced by Nagy in \cite{Nagy}, but here we produce explicit generators and that is essential to test nontriviality of K-cycles by computing the K-theory K-homology pairing. To our knowledge there are not many instances of K-theory calculations for compact quantum groups and other than the paper by Nagy there is another related work by McClanahan (\cite{McClanahan}) where he computes the K-groups of the universal $C^*$-algebra generated by a unitary matrix and shows that the associated $K_1$ is generated by the defining unitary itself. This raises the question whether something similar holds for compact matrix quantum groups, namely whether te defining unitary of a compact matrix quantum group is nontrivial in $K_1$. For quantum SU(2) this was remarked by Connes (\cite{ConSU_q2}). 
Here we not only prove that the defining unitary of  quantum SU(3) is non-trivial we show that this is part of a generating set of elements for the $K_1$.
\section{The quantum steifel manifolds and their irreducible representations }
The quantum steifel manifold 
$S_{q}^{n,m}$ was introduced by Vainerman and Podkolzin in \cite{VP}. Throughout we assume that $q \in (0,1)$.
Recall that the $C^*$ algebra $C(SU_{q}(n))$ is the universal unital $C^*$
algebra generated by $n^{2}$ elements $u_{ij}$ satisfying the following
condition 
\begin{equation}
 \sum_{k=1}^{n}u_{ik}u_{jk}^{*}=\delta_{ij}~~,~~ \sum_{k=1}^{n}u_{ki}^{*}u_{kj}
=\delta_{ij}
\end{equation}
\begin{equation}
\sum_{i_{1}=1}^{n}\sum_{i_{2}=1}^{n} \cdots \sum_{i_{n}=1}^{n}
E_{i_{1}i_{2}\cdots i_{n}}u_{j_{1}i_{1}}\cdots u_{j_{n}i_{n}} =
E_{j_{1}j_{2}\cdots j_{n}}
\end{equation}
where
\begin{displaymath}
\begin{array}{lll}
E_{i_{1}i_{2}\cdots i_{n}}&:=& \left\{\begin{array}{lll}
                                 0 &if& i_{1},i_{2},\cdots i_{n} \text{~are not
distinct} \\
                                 (-q)^{\ell(i_{1},i_{2},\cdots,i_{n})}  
                                 \end{array} \right. 
\end{array}
\end{displaymath}
where for a permutation $\sigma$ on $\{1,2,\cdots,n \}$, $\ell(\sigma)$ denotes
its length. The $C^{*}$ algebra $C(SU_{q}(n))$ has a compact quantum group
structure with the comultiplication $\Delta$ given by 
\begin{displaymath}
\Delta(u_{ij}):= \sum_{k} u_{ik}\otimes u_{kj}.
\end{displaymath} 
 
 Let $1 \leq m \leq n-1$. Call the generators of $SU_{q}(n-m)$ as $v_{ij}$. The
map $\phi:C(SU_{q}(n)) \to C(SU_{q}(n-m))$ defined by
 \begin{eqnarray}
\label{subgroup}
 \phi(u_{ij}) &:=& \left\{\begin{array}{ll}
                        v_{ij}&\text{if}~ 1\leq i,j \leq n-m ,\\
                        \delta_{ij} & \text{otherwise}.
                        \end{array} \right.
\end{eqnarray}
is a surjective unital $C^{*}$ algebra homomorphism such that $\Delta \circ \phi
= (\phi \otimes \phi)\Delta$. In this way the quantum group  $SU_{q}(n-m)$ is a
subgroup of the quantum group $SU_{q}(n)$. The $C^{*}$ algebra of the quotient
$SU_{q}(n)/SU_{q}(n-m)$ is defined as
\begin{displaymath}
C(SU_{q}(n)/SU_{q}(n-m)):=\{a \in C(SU_{q}(n)):~(\phi \otimes 1)\Delta(a)=
1\otimes a\}.
\end{displaymath} 
We refer to \cite{VP} for the proof of the following proposition.

\begin{ppsn}
The $C^{*}$ algebra $C(SU_{q}(n)/SU_{q}(n-m))$ is generated by the last $m$ rows
of the matrix $(u_{ij})$ i.e. by the set $\{u_{ij}:n-m+1 \leq i \leq n\}$.
\end{ppsn}
In \cite{VP} the quotient space $SU_{q}(n)/SU_{q}(n-m)$ is called a quantum
steiffel manifold and is denoted by $S_{q}^{n,m}$. We will also use the same
notation from now on.

The irreducible representations of the $C^{*}$
algebra $C(S_{q}^{n,m})$ was described in \cite{VP}. First we recall the
irreducible representations of $C(SU_{q}(n))$ as in \cite{So1}. The one
dimensional representations of $C(SU_{q}(n))$ are paramatrised by the  torus
$\mathbb{T}^{n-1}$. We consider $\mathbb{T}^{n-1}$ as a subset of
$\mathbb{T}^{n}$ under the inclusion $(t_{1},t_{2},\cdots t_{n-1}) \to
(t_{1},t_{2},\cdots,t_{n-1},t_{n})$ where $t_{n}:=\prod_{i=1}^{n-1}\bar{t}_{i}$.
For $t:=(t_{1},t_{2},\cdots,t_{n}) \in \mathbb{T}^{n-1}$, let
$\tau_{t}:C(SU_{q}(n)) \to \mathbb{C}$ be defined as
$\tau_{t}(u_{ij}):=t_{n-i+1}\delta_{ij}$. Then $\tau_{t}$ is a $*$ algebra
homomorphism. Moreover the set $\{\tau_{t}:t \in \mathbb{T}^{n-1}\}$ forms a
complete set of mutually inequivalent one dimensional representations of
$C(SU_{q}(n))$.

Let us denote the transposition $(i,i+1)$ by $s_{i}$. The map
$\pi_{s_{i}}:C(SU_{q}(n)) \to B(\ell^{2}(\mathbb{N}))$ defined on the generators
$u_{rs}$ as follows 
\begin{displaymath}
\begin{array}{lll}
\pi_{s_{i}}(u_{rs})&:= & \left\{ \begin{array}{lll}
                                 \sqrt{1-q^{2N+2}}S & \text{if}~~ r=i,s=i,\\
                                 -q^{N+1} & \text{if}~~ r=i,s=i+1,\\
                                  q^{N} & \text{if}~~ r=i+1,s=i,\\
                                  S^{*}\sqrt{1-q^{2N+2}} &
\text{if}~~r=i+1,s=i+1,\\
                                  \delta_{ij} & \text{otherwise}
                                \end{array} \right.
\end{array}
\end{displaymath}
is a $*$ algebra homomorphism. For any two representations $\phi$ and $\xi$ of
$C(SU_{q}(n))$, let $\phi*\xi:=(\phi\otimes \xi)\Delta$. For $\omega \in S_{n}$,
let $\omega=s_{i_{1}}s_{i_{2}}\cdots s_{i_{k}}$ be a reduced expression. Then
the representation $\pi_{\omega}:=\pi_{s_{i_{1}}}*\pi_{s_{i_{2}}}*\cdots
*\pi_{s_{i_{k}}}$ is an irreducible representation and upto unitary equivalence
the representation $\pi_{\omega}$ is independent of the reduced expression. For
$t \in \mathbb{T}^{n-1}$ and $\omega \in S_{n}$ let $\pi_{t,\omega}:=
\tau_{t}*\pi_{\omega}$. We refer to \cite{So1} for the proof of the following
theorem.

\begin{thm}
The set $\{\pi_{t,\omega}: t \in \mathbb{T}^{n-1},\omega \in S_{n}\}$ forms a
complete set of mutually inequivalent irreducible representations of
$C(SU_{q}(n))$.
\end{thm}
In \cite{VP} the irreducible representations of $C(S_{q}^{n,m})$ are studied and
we recall them here. We embedd $\mathbb{T}^{m}$ into $\mathbb{T}^{n-1}$ via the
map $t=(t_{1},t_{2},\cdots,t_{m}) \to
(t_{1},t_{2},\cdots,t_{m},1,1,\cdots,1,t_{n})$ where
$t_{n}:=\prod_{i=1}^{m}\bar{t}_{i}$. For a permutation $\omega \in S_{n}$, let
$\omega^{s}$ be the permutation in the coset $S_{n-m}\omega$ with the least
possible length. We denote the restriction of the representation
$\pi_{t,\omega}$  to the subalgbera $C(S_{q}^{n,m})$ by $\pi_{t,\omega}$ itself.
Then we have the following theorem whose proof can be found in \cite{VP}
\begin{thm}
The set $\{\pi_{t,\omega^{s}}:t \in \mathbb{T}^{m},\omega \in S_{n}\}$ forms a
complete set of mutually inequivalent irreducible representations of
$C(S_{q}^{n,m})$.
\end{thm} 
\section{Composition sequences}
In this section we derive certain exact sequences analogous to that of Theorem 4
in \cite{Sh1}. We then apply the six term sequence in K theory to compute the K
groups of $C(S_{q}^{n,2})$. 
\begin{lmma}
Let $t \in \mathbb{T}^{m}$ and $\omega:=s_{n-1}s_{n-2}\cdots s_{n-k}$. Then 
Image of $ C(S_{q}^{n,m})$ under the homomorphism $\pi_{t,\omega}$  contains the algebra of compact operators
$\mathcal{K}(\ell^{2}(\mathbb{N}^{k}))$.
\end{lmma}
\textit{Proof.} 
Since
$\pi_{t,\omega}(C(S_{q}^{n,m}))=\pi_{\omega}(C(S_{q}^{n,m}))$, it is enough to
show that $\mathcal{K}(\ell^{2}(\mathbb{N}^{k})) \subset
\pi_{\omega}(C(S_{q}^{n,m}))$. We prove this result by induction on $n$. Since
$\pi_{\omega}(u_{nn}):=S^{*}\sqrt{1-q^{2N+2}}\otimes 1$, it follows that
$S\otimes 1 \in \pi_{\omega}(C(S_{q}^{n,m})$. Hence
$\mathcal{K}(\ell^{2}(\mathbb{N}))\otimes 1 \subset \pi_{w}(C(S_{q}^{n,m})$.
Thus the result is true if $n=2$. Next observe that for $1 \leq i \leq n-1$,
$(p\otimes 1)\pi_{\omega}(u_{n,i}):= p \otimes \pi_{\omega^{'}}(v_{n-1,i})$ 
where $\omega^{'}:=s_{n-2}s_{n-3}\cdots s_{n-k}$ and $(v_{ij})$ denotes the
generators of $C(SU_{q}(n-1))$. Hence $\pi_{\omega}(C(S_{q}^{n,m}))$ contains
the algebra $p\otimes \pi_{\omega^{'}}(C(S_{q}^{n-1,m}))$. Now  by induction
hypothesis, it follows that $\pi_{\omega}(C(S_{q}^{n,m}))$ contains $p \otimes
\mathcal{K}(\ell^{2}(\mathbb{N}^{k-1}))$. Since $\pi_{\omega}(C(S_{q}^{n,m}))$
contains $\mathcal{K}(\ell^{2}(\mathbb{N})) \otimes 1$ and $p \otimes
\mathcal{K}(\ell^{2}(\mathbb{N}^{k-1}))$, it follows that
$\pi_{\omega}(C(S_{q}^{n,m}))$ contains the algebra of compact operators. This
completes the proof. \hfill$\Box$

Let $w$ be a word on $s_{1},s_{2},\cdots s_{n}$ say $w:=s_{i_{1}}s_{i_{2}}\cdots
s_{i_{n}}$ (not necessarily a reduced expression). Define
$\psi_{w}:=\pi_{s_{i_{1}}}*\pi_{s_{i_{2}}}*\cdots \pi_{s_{i_{r}}}$ and for $t\in
\mathbb{T}^{n}$, let $\psi_{t,w}:= \tau_{t}*\psi_{w}$. Observe that the image of
 $\psi_{t,w}$ is contained in $\tau^{\otimes r}$.  We prove that if $w^{'}$ is a
'subword' of $w$ then $\psi_{t,w^{'}}$ factors through $\psi_{t,w}$.
\begin{ppsn}
\label{factorisation}
Let $w=w_{1}s_{k}w_{2}$ be a word on $s_{1},s_{2},\cdots s_{n}$. Denote the word
$w_{1}w_{2}$ by $w^{'}$. Let $t \in \mathbb{T}^{m}$ be given. Then there exists
a * homomorphism $\epsilon:\psi_{t,w}(C(S_{q}^{n,m})) \to
\psi_{t,w^{'}}(C(S_{q}^{n,m}))$ such that $\psi_{t,w^{'}}= \epsilon \circ
\psi_{t,w}$.
\end{ppsn}
\textit{Proof.} For a word $u$ on $s_{1},s_{2},\cdots,s_{n}$, let $\ell(u)$
denote it's length. Then $\psi_{t,w}(C(S_{q}^{n,m}))$ is contained in
$\tau^{\otimes \ell(w_{1})}\otimes \tau \otimes \tau^{\otimes \ell(w_{2})}$. Let
$\epsilon$ denote the restriction of $1\otimes \sigma \otimes 1$ to
$\psi_{t,w}(C(S_{q}^{n,m}))$ where $\sigma:\tau \to \mathbb{C}$ is the
homomorphism for which $\sigma(S)=1$.
\begin{displaymath}
\psi_{t,w}(u_{rs})= \sum_{j_{1},j_{2}}\psi_{t,w_{1}}(u_{rj_{1}})\otimes
\pi_{s_{k}}(u_{j_{1}j_{2}})\otimes \psi_{w_{2}}(u_{j_{2}s}).
\end{displaymath}
Since $\sigma(\pi_{s_{k}}(u_{j_{1}j_{2}}))=\delta_{j_{1}j_{2}}$, it follows
that 
\begin{displaymath}
\epsilon \circ \psi_{t,w}(u_{rs})= \sum_{j} \psi_{t,w_{1}}(u_{rj})\otimes
\psi_{w_{2}}(u_{js})=\psi_{t,w^{'}}(u_{rs}).
\end{displaymath} 
This completes the proof. \hfill $\Box$

Let $w$ be a word on $s_{1},s_{2},\cdots s_{n}$. Then for $n-m+1 \leq i \leq n$
and $1 \leq j \leq n$, the map $\mathbb{T}^{m}:t \to \psi_{t,w}(u_{ij}) \in
\tau^{\otimes \ell(w)}$ is continuous. Thus we get a homomorphism
$\chi_{w}:C(S_{q}^{n,m}) \to C(\mathbb{T}^{m})\otimes \tau^{\otimes \ell(w)}$
such that $\chi_{w}(a)(t)=\psi_{t,w}(a)$ for all $a \in C(S_{q}^{n,m})$.
\begin{rmrk}
\label{factorisation1}
 Clearly for a word $w$ on $s_{1},s_{2},\cdots s_{n}$ the representations
$\psi_{t,w}$ factors through $\chi_{w}$. One can also prove as in lemma
\ref{factorisation} that if $w^{'}$ is a 'subword' of $w$ then $\chi_{w^{'}}$
factors through $\chi_{w}$. 
\end{rmrk}
Let us introduce some notations. Denote the permutation $s_{j}s_{j-1}\cdots
s_{i}$ by $\omega_{j,i}$ for $j\geq i$. If $j>i$ we let $\omega_{j,i}:=1$. For
$1 \leq k \leq n$, let $\omega_{k}:=\omega_{n-m,1}\omega_{n-m+1,1}\cdots
\omega_{n-1,n-k+1}$.
\begin{thm}
The homomorphism  $\chi_{\omega_{n}}:C(S_{q}^{n,m})\to C(\mathbb{T}^{m})\otimes
\tau^{\otimes \ell(\omega_{n})}$ is faithful.
\end{thm}
\textit{Proof.} If $\omega_{0} \in S_{n}$ then $\omega_{0}^{s}$( the
representative in $S_{n-m}\omega_{0}$ with the shortest length) is a 'subword'
of $\omega_{n}$. Hence by remark \ref{factorisation1} every irreducible
representation of $C(S_{q}^{n,m})$ factors through $\chi_{\omega_{n}}$. Thus it
follows that $\chi_{\omega_{n}}$ is faithful. This completes the proof. \hfill
$\Box$

For $1\leq k \leq n$, Let $C(S_{q}^{n,m,k}):=\chi_{\omega_{k}}(C(S_{q}^{n,m}))$.
Then $C(S_{q}^{n,m,k}) \subset C(S_{q}^{n,m,1})\otimes \tau^{\otimes(k-1)}$. For
$2 \leq k \leq n$, let $\sigma_{k}$ denote the restriction of $(1\otimes
1^{\otimes(k-2)}\otimes \sigma)$ to $C(S_{q}^{n,m,k})$. Then the image of
$\sigma_{k}$ is $C(S_{q}^{n,m,k-1})$. We determine the kernel of $\sigma_{k}$ in
the next proposition. We need the following two lemmas.
\begin{lmma}
\label{ideal of oddspheres}
The algebra $\chi_{\omega_{n-1,n-k}}(C(S_{q}^{n,1}))$ contains
$C^{*}(t_{1})\otimes \mathcal{K}(\ell^{2}(\mathbb{N}^{k}))$ which is isomorphic
to $C(\mathbb{T})\otimes \mathcal{K}(\ell^{2}(\mathbb{N}^{k}))$.
\end{lmma}
\textit{Proof.} Note that $\chi_{\omega_{n-1,n-k}}(u_{nn})= t_{1}\otimes
S^{*}\sqrt{1-q^{2N+2}}\otimes 1$. Hence it follows that the operator  $1\otimes
\sqrt{1-q^{2N+2}} \otimes 1=\chi_{\omega_{n-1,n-k}}(u_{nn}^{*}u_{nn})$ lies in
the algebra $\chi_{\omega_{n-1,n-k}}(C(S_{q}^{n,1}))$. As $\sqrt{1-q^{2N+2}}$ is
invertible, one has $t_{1}\otimes S^{*}\otimes 1 \in
\chi_{\omega_{n-1,n-k}}(C(S_{q}^{n,1}))$. Thus the projection $1\otimes p
\otimes 1$ is in the algebra $C(S_{q}^{n,1,k+1})$. Now observe that for $1 \leq
s \leq n-1$, one has 
\begin{equation}
\label{eq}
(1\otimes p \otimes 1)\chi_{\omega_{n-1,n-k}}(u_{ns})= t_{1}\otimes p \otimes
\pi_{\omega_{n-2,n-k}}(v_{n-1,s})
\end{equation} where $(v_{ij})$ are the generators of $C(SU_{q}(n-1))$. If $n=2$
then $k=1$ and what we have shown is that $C(S_{q}^{2,1,2})$ contains
$t_{1}\otimes S^{*}$ and $t_{1}\otimes p$. Hence one has $C^{*}(t_{1})\otimes
\mathcal{K}$ is contained in the algebra $C(S_{q}^{2,1,2})$.

Now we can complete the proof by induction on $n$. Equation \ref{eq} shows that
$C^{*}(t_{1})\otimes p \otimes \mathcal{K}^{\otimes(k-1)}$ is contained in the
algebra $C(S_{q}^{n,1,k+1})$ and we also have $t_{1}\otimes S^{*}\otimes 1 \in
C(S_{q}^{n,1,k+1})$. Hence it follows that $C^{*}(t_{1})\otimes
\mathcal{K}^{\otimes k}$ is contained in the algebra $C(S_{q}^{n,1,k+1})$. This
completes the proof. \hfill $\Box$

\begin{lmma}
\label{killing}
 Given $1\leq s \leq n$, there exists compact operators $x_{s},y_{s}$ such that
$x_{s}\pi_{\omega_{n-1,n-k}}(u_{js})y_{s}= \delta_{js}(p\otimes p \otimes \cdots
\otimes p)$ where $p:=1-S^{*}S$.
\end{lmma}
\textit{Proof.} Let $1\leq s\leq n$ be given. Note that the operator
$\omega_{n-1,n-k}(u_{ss})=z_{1}\otimes z_{2} \otimes \cdots z_{k}$ where $z_{i}
\in \{1,\sqrt{1-q^{2N+2}}S,S^{*}\sqrt{1-q^{2N+2}}\}$. Define $x_{i},y_{i}$ as
follows
\begin{displaymath}
\begin{array}{lll}
 x_{i}:=& \left\{ \begin{array}{lll}
                 p & \mbox{if} & z_{i}=1, \\
                 p & \mbox{if} & z_{i}=\sqrt{1-q^{2N+2}}S ,\\
                 (1-q^{2})^{-\frac{1}{2}}pS & \mbox{if} &
z_{i}=S^{*}\sqrt{1-q^{2N+2}}.
\end{array} \right.
\end{array}
\end{displaymath}
\begin{displaymath}
\begin{array}{lll}
 y_{i}:=& \left\{ \begin{array}{lll}
                 p & \mbox{if} & z_{i}=1, \\
                 (1-q^{2})^{-\frac{1}{2}}S^{*}p & \mbox{if} &
z_{i}=\sqrt{1-q^{2N+2}}S ,\\
                 p & \mbox{if} & z_{i}=S^{*}\sqrt{1-q^{2N+2}}.
\end{array} \right.
\end{array}
\end{displaymath}
Then $x_{i}z_{i}y_{i}=p$ for $1\leq i \leq k$. Now let $x_{s}:=x_{1}\otimes
x_{2}\otimes \cdots x_{k}$ and $y_{s}:=y_{1}\otimes y_{2} \otimes \cdots y_{k}$.
Then $x_{s}\chi_{\omega_{n-1,n-k}}(u_{ss})=\underbrace{p \otimes p\otimes \cdots
\otimes p}_{k times}$. Let $j\neq s$ be given. Then
$\chi_{\omega_{n-1,n-k}}(u_{js})=a_{1}\otimes a_{2} \otimes \cdots \otimes
a_{k}$ where $a_{i} \in
\{1,\sqrt{1-q^{2N+2}}S,S^{*}\sqrt{1-q^{2N+2}},-q^{N+1},q^{N} \}$. Since $j\neq
s$, there exists an $i$ such that $a_{i} \in \{q^{N},-q^{N+1}\}$. Let $r$ be the
largest integer for which $a_{r} \in \{q^{N},-q^{N+1}\}$. Then $z_{r}\neq 1$.
Hence $x_{r}a_{r}y_{r}=0$. Thus $x_{s}\chi_{\omega_{n-1,n-k}}(u_{js})y_{s}=0$.
This completes the proof. \hfill $\Box$
\begin{ppsn}
\label{sheu}
Let $2 \leq k \leq n$. Then $C(S_{q}^{n,m,1})\otimes
\mathcal{K}(\ell^{2}(\mathbb{N}))^{\otimes(k-1)}$ is contained in the algebra
$C(S_{q}^{n,m,k})$. Moreover the kernel of the homomorphism $\sigma_{k}$ is
exactly $C(S_{q}^{n,m,1})\otimes
\mathcal{K}(\ell^{2}(\mathbb{N}))^{\otimes(k-1)}$. Thus we have the exact
sequence 
\begin{displaymath}
0 \longrightarrow C(S_{q}^{n,m,1})\otimes \mathcal{K}^{\otimes(k-1)}
\longrightarrow C(S_{q}^{n,m,k}) \stackrel{\sigma_{k}} \longrightarrow
C(S_{q}^{n,m,k-1}) \longrightarrow 0.
\end{displaymath}
\end{ppsn}
\textit{Proof.} First we prove that $C(S_{q}^{n,m,1})\otimes
\mathcal{K}^{\otimes(k-1)}$ is contained in the algebra $C(S_{q}^{n,m,k})$. For
$a \in C(S_{q}^{n,1})$ one has $\chi_{\omega_{k}}(a):=1 \otimes
\chi_{\omega_{n-1,n-k+1}}(a)$, it follows from lemma \ref{ideal of oddspheres}
that $C(S_{q}^{n,m,k})$ contains $1 \otimes
\mathcal{K}(\ell^{2}(\mathbb{N}^{k-1}))$. Let $n-m+1 \leq r \leq m$ and $1 \leq
s \leq n$ be given. Then note that 
\begin{displaymath}
\chi_{\omega_{k}}(u_{rs})=\sum_{j=1}^{n}\chi_{\omega_{1}}(u_{rj})\otimes
\pi_{\omega_{n-1,n-k+1}}(u_{js}).
\end{displaymath}
Hence by lemma \ref{killing}, there exists $x_{s},y_{s} \in C(S_{q}^{n,m,k})$
such that $x_{s}\chi_{\omega_{k}}(u_{rs})y_{s}:=
\chi_{\omega_{1}}(u_{rs})\otimes p^{\otimes(k-1)}$ where
$p^{\otimes(k-1)}:=p\otimes p\otimes \cdots \otimes p$. Thus we have shown that
$C(S_{q}^{n,m,k})$ contains $1\otimes \mathcal{K}^{\otimes(k-1)}$ and
$C(S_{q}^{n,m,1})\otimes p^{\otimes(k-1)}$. Hence $C(S_{q}^{n,m,k})$ contains
$C(S_{q}^{n,m,1})\otimes \mathcal{K}^{\otimes(k-1)}$.

Clearly $\sigma_{k}$ vanishes on $C(S_{q}^{n,m,1})\otimes
\mathcal{K}^{\otimes(k-1)}$. Let $\pi$ be an irreducible representation of
$C(S_{q}^{n,m,k})$ which vanishes on the ideal $C(S_{q}^{n,m,1})\otimes
\mathcal{K}^{\otimes(k-1)}$. Then $\pi\circ \chi_{\omega_{k}}$ is an irreducible
representation of $C(S_{q}^{n,m})$. Hence $\pi \circ
\chi_{\omega_{k}}=\pi_{t,\omega}$ for some $\omega$ of the form
$\omega_{n-m,i_{1}}\omega_{n-m+1,i_{2}}\cdots \omega_{n-1,i_{n-m}}$ and $t \in
\mathbb{T}^{m}$. Since $\pi \circ \chi_{\omega_{k}}(u_{n,n-k+1})=0$, it follows
that $\pi_{t,w}(u_{n,n-k+1})= 0$. But one has
$\pi_{t,\omega}(u_{n,n-k+1})=t_{n}(1\otimes
\pi_{\omega_{n-1,i_{n-m}}}(u_{n,n-k+1}))$. Hence $i_{n-m}>n-k+1$. In other words
$\omega$ is a subword of $\omega_{k-1}$. Thus $\pi \circ \chi_{\omega_{k}}$
factors through $\chi_{\omega_{k-1}}$. In other words there exists a
representation $\rho$ of $C(S_{q}^{n,m,k-1})$ such that $\pi \circ
\chi_{\omega_{k}}=\rho \circ \chi_{\omega_{k-1}}$. Since
$\chi_{\omega_{k-1}}=\sigma_{k}\circ \chi_{\omega_{k}}$, it follows that
$\pi=\rho\circ \sigma_{k}$. Thus we have shown that every irreducible
representation of $C(S_{q}^{n,m,k})$ which vanishes on the ideal
$C(S_{q}^{n,m,1})\otimes \mathcal{K}^{\otimes(k-1)}$ factors through
$\sigma_{k}$. Hence the kernel of $\sigma_{k}$ is exactly the ideal
$C(S_{q}^{n,m,1})\otimes \mathcal{K}^{\otimes(k-1)}$. This completes the proof.
\hfill$\Box$

We apply the six term exact sequence in K theory to the exact sequence in
proposition \ref{sheu} to compute the $K$ groups of $C(S_{q}^{n,2,k})$ for $1
\leq k \leq n$. In the next section we breifly recall the product operation in
$K$ theory.

\section{The operation P}
 Let $A$ and $B$ be $C^{*}$ algebras. Then we have the following product maps.
\begin{align*}
  K_{0}(A) \otimes K_{0}(B) & \to K_{0}(A\otimes B),\\
  K_{1}(A) \otimes K_{0}(B) & \to K_{1}(A \otimes B), \\
  K_{0}(A) \otimes K_{1}(B) & \to K_{1}(A \otimes B) ,\\
  K_{1}(A) \otimes K_{1}(B) & \to K_{0}(A \otimes B). \\
\end{align*}
 The first map is defined as $[p] \otimes [q] \to [p \otimes q]$. The second one
is defined as $[u] \otimes [p] \to [u \otimes p + 1-1 \otimes p]$. The third map
is defined in the same manner and the fourth one is defined using Bott
periodicity and using the first product.
In fact we have the following formula for the last product. We refer to Appendix
of \cite{Con}.

Let $h: \mathbb{T}^{2} \to P_{1}(\mathbb{C}):= \{p \in Proj(M_{2}(\mathbb{C}):
trace(p)=1 \}$ be a degree one map. Then given unitaries $u \in M_{p}(A)$ and $v
\in M_{q}(B)$ the product $[u] \otimes [v]$ is given by $[h(u,v)]-[e_{0}]$ where
$\begin{array}{ll}
  e_{0}=& \left ( \begin{array}{ll}
                   1 & 0 \\
                   0 & 0 
                 \end{array} \right )
\end{array} \in M_{2}(M_{pq}(A \otimes B))$ and $h(u,v)$ is the matrix with
entries $h_{ij}(u \otimes 1, 1 \otimes v)$.
We denote the image of $[x] \otimes [y]$ by $[x] \otimes [y]$ itself. 
Now let $A$ be a unital commutative $C^*$ algebra. Then the multiplication $m:A
\otimes A \to A$ is a $C^{*}$ algebra homomorphism. Hence we get a map at the K
theory level from $K_{1}(A) \otimes K_{1}(A) \to K_{0}(A)$. 

Suppose $U$ and $V$ are two commuting unitaries in a $C^{*}$ algebra $B$. Let
$A:=C^{*}(U,V)$. Then $A$ is commutative. Define
\begin{equation*}
P(U,V): = K_{0}(m)([U] \otimes [V])
\end{equation*}
 which is an element in $K_{0}(A)$ which we can think of as an element in
$K_{0}(B)$ by composing with the inclusion map. From the formula that we just
recalled from \cite{Con} the following properties are clear :
\begin{enumerate}
\item If $U$ and $V$ are commuting unitaries in $A$ and $p$ is a rank one
projection in $\mathcal{K}$ we have $P(U \otimes p + 1-1 \otimes p, V \otimes p
+ 1-1\otimes p):= P(U,V) \otimes p$.
\item If $U$ and $V$ are commuting unitaries and $p$ is a projection that
commutes with $U$ and $V$ then $P(U,Vp+1-p)=P(Up+1-p,Vp+1-p)$.
\item If $\phi:A \to B$ is a unital homomorphism and if $U$ and $V$ are
commuting unitaries in $A$ then $K_{0}(\phi)(P(U,V))= P(\phi(U),\phi(V))$.
\item If $U$ is a unitary in $A$ then $P(U,U)=0$. For $P_{1}(\mathbb{C})$ is
simply connected, it follows that the matrix $h(U,U)$ is path connected to a
rank one projection in $M_{2}(\mathbb{C})$. Hence $P(U,U)=0$.
\end{enumerate} 
We need the following lemma in the six term computation. Let $z_{1}\otimes 1$
and $1\otimes z_{2}$ be the generating unitaries of $C(\mathbb{T})\otimes
C(\mathbb{T})$. Then $K_{0}(C(\mathbb{T}^{2}))$ is isomorphic to
$\mathbb{Z}^{2}$ and is generated by $1,P(z_{1}\otimes 1,1\otimes z_{2})$.
\begin{lmma}
\label{index}
Consider the exact sequence 
\begin{displaymath}
0 \longrightarrow C(\mathbb{T})\otimes \mathcal{K} \longrightarrow
C(\mathbb{T})\otimes \tau \longrightarrow C(\mathbb{T})\otimes C(\mathbb{T})
\longrightarrow 0
\end{displaymath}
and the six term sequence in $K$ theory.
\begin{equation*}
\def\labelstyle{\scriptstyle}
\xymatrix@C=25pt@R=20pt{
K_0(C(\mathbb{T})\otimes \mathcal{K})\ar[r]& K_0(C(\mathbb{T})\otimes
\tau)\ar[r]& K_0(C(\mathbb{T})\otimes C(\mathbb{T})) \ar[d]_{\delta} \\
K_1(C(\mathbb{T})\otimes C(\mathbb{T}))\ar[u]^{\partial}&
K_1(C(\mathbb{T})\otimes \tau)\ar[l] & K_1(C(\mathbb{T})\otimes
\mathcal{K})\ar[l] }
\end{equation*}
Then the subgroup generated by $\delta(P(z_{1}\otimes 1,1\otimes z_{2}))$
coincides with the group generated by $z_{1}\otimes p + 1-1\otimes p$ which is
$K_{1}(C(\mathbb{T})\otimes \mathcal{K}) \cong \mathbb{Z}$.
\end{lmma}
\textit{Proof.} The toeplitz map $\epsilon:\tau \to C(\mathbb{T})$ induces
isomorphism at the $K_{0}$  level. Thus by Kunneth theorem, it follows that the
image of $K_{0}(1\otimes \epsilon)$ is $K_{0}(C(\mathbb{T}))\otimes
K_{0}(C(\mathbb{T}))$ which is the subgroup generated by $[1]$. Now the
inclusion $0 \to \mathcal{K} \to \tau$ induces the zero map at the $K_{0}$ level
and hence again by Kunneth theorem the inclusion $0 \to C(\mathbb{T})\otimes
\mathcal{K} \to C(\mathbb{T})\otimes \tau$ induces zero map at the $K_{1}$
level. Hence the image of $\delta$ is $K_{1}(C(\mathbb{T})\otimes \tau)$. This
completes the proof.
\begin{crl}
\label{toeplitz trick}
Let 
\begin{displaymath}
0 \longrightarrow I \longrightarrow A  \stackrel{\phi} \longrightarrow B
\longrightarrow 0
\end{displaymath}
be a short exact sequence of $C^{*}$ algebras. Consider the six term sequence in
$K$ theory.
\begin{equation*}
\label{indexmap}
\def\labelstyle{\scriptstyle}
\xymatrix@C=25pt@R=20pt{
K_0(I)\ar[r]& K_0(A)\ar[r]^{K_0(\phi)}& K_0(B) \ar[d]_{\delta} \\
K_1(B)\ar[u]^{\partial}& K_1(A)\ar[l]^{K_{1}(\phi)}& K_1(I)\ar[l] 
}
\end{equation*}
Suppose that $U$ and $V$ are two commuting unitaries in $B$ such that there
exists a unitary $X$ and an isometry $Y$ such that $\phi(X)=U$ and $\phi(Y)=V$.
Also assume that $X$ and $Y$ commute. Then the subgroup generated by
$\delta(P(U,V))$ coincides with the subgroup generated by the unitary
$X(1-YY^{*})+ YY^{*}$ in $K_{1}(I)$.
\end{crl}
\textit{Proof.} Since $C(\mathbb{T})$ is the universal $C^{*}$ algebra generated
by a unitary and $\tau$ is the universal $C^{*}$ algebra generated by an
isometry, there exists homomorphisms $\Phi:C(\mathbb{T})\otimes \tau \to A$ and
$\Psi:C(\mathbb{T})\otimes C(\mathbb{T}) \to B$ such that 
\begin{displaymath}
\begin{array}{ll}
\Phi(z_{1}\otimes 1):= X ,\\
\Phi(1\otimes S^{*}):= Y ,\\
\Psi(z_{1} \otimes 1):= U, \\
\Psi(1\otimes z_{2}):= V.
\end{array}
\end{displaymath}
Hence we have the following commutative diagram.
\begin{equation*}
\label{indexmap}
\def\labelstyle{\scriptstyle}
\xymatrix@C=25pt@R=20pt{
0\ar[r] & C(\mathbb{T})\otimes \mathcal{K}\ar[d]^{\Phi}\ar[r] &
C(\mathbb{T})\otimes \tau \ar[d]^{\Phi}\ar[r] & C(\mathbb{T})\otimes
C(\mathbb{T})\ar[d]^{\Psi}\ar[r] & 0\\
0 \ar[r] & I \ar[r] & A \ar[r]^{\phi} & B \ar[r] & 0  
}
\end{equation*}
Now by the functoriality of $\delta$ and $P$, it follows that
$\delta(P(U,V))=K_{1}(\Phi)(\delta(P(z_{1}\otimes 1,1\otimes z_{2})))$. Hence by
lemma \ref{index}, it follows that the subgroup generated by $\delta(P(U,V))$ is
the subgroup generated by $\Phi(z_{1}\otimes p+1-1\otimes p)$ in $K_{1}(I)$.
Note that $\Phi(z_{1}\otimes p+1-1\otimes p)=X(1-YY^{*})+YY^{*}$. This completes
the proof. \hfill $\Box$.
\section{$K$ groups of $C(S_{q}^{n,2,k})$ for $k<n$}
In this section we compute the $K$ groups of $C(S_{q}^{n,2,k})$ for $1\leq k <n$
by applying the six term sequence in $K$ theory to the exact sequence in
\ref{sheu}. Let us fix some notations. If $q$ is a projection in
$\ell^{2}(\mathbb{N})$ then $q_{r}$ denotes the projection $\underbrace{q\otimes
q\otimes \cdots q}_{r~times}$ in $\ell^{2}(\mathbb{N}^{r})$. Let us define the
unitaries $U_{k},V_{k},u_{k},v_{k}$ as follows.
\begin{displaymath}
\begin{array}{ll}
U_{k}:=&t_{1}\otimes 1_{n-2}\otimes p_{k-1}+ 1-1\otimes 1_{n-2}\otimes p_{k-1},\\
V_{k}:=&t_{2}\otimes p_{n-2}\otimes 1_{k-1}+ 1-1\otimes p_{n-2} \otimes
1_{k-1},\\
u_{k}:=&t_{1}\otimes p_{n-2}\otimes p_{k-1}+ 1-1\otimes p_{n-2}\otimes p_{k-1},\\
v_{k}:=&t_{2}\otimes p_{n-2}\otimes p_{k-1}+ 1-1\otimes p_{n-2}\otimes p_{k-1}.\\
\end{array}
\end{displaymath}
First let us note that the operators $U_{k},V_{k},u_{k},v_{k}$ lies in the
algebra $C(S_{q}^{n,2,k})$. For,
\begin{displaymath}
\begin{array}{ll}
U_{k}=&1_{\{1\}}(u_{n,n-k+1}u_{n,n-k+1}^{*})u_{n,n-k+1}+
1-1_{\{1\}}(u_{n,n-k+1}u_{n,n-k+1}^{*}),\\
V_{k}=&1_{\{1\}}(u_{n-1,1}u_{n-1,1}^{*})u_{n-1,1}+
1-1_{\{1\}}(u_{n-1,1}u_{n-1,1}^{*}),\\
u_{k}=&1_{\{1\}}(u_{n,n-k+1}u_{n,n-k+1}^{*}u_{n-1,1}u_{n-1,1}^{*})u_{n,n-k+1}+ 1
-1_{\{1\}}(u_{n,n-k+1}u_{n,n-k+1}^{*}u_{n-1,1}u_{n-1,1}^{*}),\\
v_{k}=&1_{\{1\}}(u_{n,n-k+1}u_{n,n-k+1}^{*}u_{n-1,1}u_{n-1,1}^{*})u_{n-1,1}+1-1_
{\{1\}}(u_{n,n-k+1}u_{n,n-k+1}^{*}u_{n-1,1}u_{n-1,1}^{*}).
\end{array}
\end{displaymath}
Note that the unitaries $U_{n},u_{n},v_{n}$ lies in the algebra
$C(S_{q}^{n,2,n})$. We start with the computation of the $K$ groups of
$C(S_{q}^{n,2,1})$.
\begin{lmma}
\label{K groups}
The $K$ groups $K_{0}(C(S_{q}^{n,2,1}))$ and $K_{1}(C(S_{q}^{n,2,1}))$ are both
isomorphic to $\mathbb{Z}^{2}$. In fact, $[U_{1}]$ and $[V_{1}]$ form a
$\mathbb{Z}$ basis for $K_{1}(C(S_{q}^{n,2,1}))$ and $[1]$ and $P(u_{1},v_{1})$
form a $\mathbb{Z}$ basis for $K_{0}(C(S_{q}^{n,2,1}))$.
\end{lmma}
\textit{Proof.} First note that $C(S_{q}^{n,2,1})$ is generated by $t_{1}\otimes
1_{n-2}$ and $t_{2}\otimes \pi_{\omega_{n-2,1}}(u_{n-1,j})$ where $1\leq j \leq
n-1$. But the $C^{*}$ algebra generated by $\{t_{2}\otimes
\pi_{\omega_{n-2,1}}(u_{n-1,j}): ~1\leq j \leq n-1\}$ is isomorphic to
$C(S_{q}^{2n-3})$. Hence $C(S_{q}^{n,2,1})$ is isomorphic to
$C(\mathbb{T})\otimes C(S_{q}^{2n-3})$. Also $K_{0}(C(S_{q}^{2n-3})$ and
$K_{1}(C(S_{q}^{2n-3}))$ are both isomorphic to $\mathbb{Z}$ with $[1]$
generating $K_{0}(C(S_{q}^{2n-3}))$ and $[t_{2}\otimes p_{n-2}+1-1\otimes
p_{n-2}]$ generating $K_{1}(C(S_{q}^{2n-3}))$.

Now by the Kunneth theorem for tensor product of $C^{*}$ algebras(See
\cite{Bla}), it follows that $C(S_{q}^{n,2,1})$ has both $K_{1}$ and $K_{0}$
isomorphic to $\mathbb{Z}^{2}$ with $[U_{1}]$ and $[V_{1}]$ generating
$K_{1}(C(S_{q}^{n,2,1}))$ and $[1]$ and $P(t_{1}\otimes 1_{n-2},t_{2}\otimes
p_{n-2}+1-1\otimes p_{n-2})$ generating $K_{0}(C(S_{q}^{n,2,1}))$. Note that the
projection $1\otimes p_{n-2}=
1_{\{1\}}(\chi_{\omega_{n-2,1}}(u_{n-1,1}u_{n-1,1}^{*}))$ is  in
$C(S_{q}^{n,2,1})$ and commutes with the unitaries $t_{1}\otimes 1_{n-2}$ and
$t_{2}\otimes p_{n-2}+1-1\otimes p_{n-2}$. Hence 
\begin{displaymath}
P(t_{1}\otimes 1_{n-2},t_{2}\otimes p_{n-2}+1-1\otimes p_{n-2})= P(u_{1},v_{1}).
\end{displaymath}
This completes the proof. \hfill $\Box$

\begin{ppsn}
Let $1\leq k <n$ be given. Then the $K_{0}(C(S_{q}^{n,2,k}))$ and
$K_{1}(C(S_{q}^{n,2,k}))$ are both isomorphic to $\mathbb{Z}^{2}$. In
particular, $[U_{k}]$ and $[V_{k}]$ form a $\mathbb{Z}$ basis for
$K_{1}(C(S_{q}^{n,2,k}))$ and $[1]$ and $P(u_{k},v_{k})$ form a $\mathbb{Z}$
basis for $K_{0}(C(S_{q}^{n,2,k}))$.
\end{ppsn}
\textit{Proof.} We prove this result by induction on $k$. The case $k=1$ is just
lemma \ref{K groups}. Now assume the result to be true for $k$. From proposition
\ref{sheu} we have the short exact sequence
\begin{displaymath}
0 \longrightarrow C(S_{q}^{n,2,1})\otimes \mathcal{K}^{\otimes(k)}
\longrightarrow C(S_{q}^{n,2,k+1}) \stackrel{\sigma_{k+1}} \longrightarrow
C(S_{q}^{n,2,k}) \longrightarrow 0
\end{displaymath}
which gives rise to the following six term sequence in $K$ theory.
\begin{equation*}
\def\labelstyle{\scriptstyle}
\xymatrix@C=25pt@R=20pt{
K_0(C(S_{q}^{n,2,1})\otimes \mathcal{K}^{\otimes k}) \ar[r]&
K_0(C(S_{q}^{n,2,k+1}))\ar[r]^{K_0(\sigma_{k+1})}& K_0(C(S_{q}^{n,2,k}))
\ar[d]_{\delta} \\
K_1(C(S_{q}^{n,2,k}) \ar[u]^{\partial} &
K_1(C(S_{q}^{n,2,k+1}))\ar[l]^{K_{1}(\sigma_{k+1})}& K_1(C(S_{q}^{n,2,1})\otimes
\mathcal{K}^{\otimes k})\ar[l] 
}
\end{equation*}
We determine $\delta$ and $\partial$ to compute the six term sequence. As
$\sigma_{k+1}(V_{k+1})=V_{k}$, it follows that $\partial([V_{k}])=0$. Since
$C(S_{q}^{n,2,k+1})$ contains the algebra $C(S_{q}^{n,2,1})\otimes
\mathcal{K}^{\otimes k}$, it follows that the operator $\tilde{X}:=t_{1}\otimes
1_{n-2}\otimes \underbrace{q^{N}\otimes q^{N}\otimes \cdots q^{N}}_{(k-1) times}
\otimes S^{*}$ is in the algebra $C(S_{q}^{n,2,1})$ as the difference
$X-\chi_{\omega_{k+1}}(u_{n,n-k+1})$ lies in the ideal $C(S_{q}^{n,2,1})\otimes
\mathcal{K}^{\otimes k}$. Let $X:=1_{\{1\}}(\tilde{X}^{*}\tilde{X})\tilde{X} +
1-1_{\{1\}}(\tilde{X}^{*}\tilde{X})$. Then $X$ is an isometry such that
$\sigma_{k+1}(X)=U_{k}$. Hence $\partial([U_{k}])=[1-X^{*}X]-[1-XX^{*}]$. Thus
$\partial([U_{k}])=-[1\otimes 1_{n-2}\otimes p_{k}]$. Thus the image of
$\partial$ is the subgroup of $K_{0}(C(S_{q}^{n,2,1})\otimes
\mathcal{K}^{\otimes k})$ generated by $[1\otimes 1_{n-2}\otimes  p_{k}]$ and
the kernel is $[V_{k}]$.

Next we compute $\delta$. Since $\sigma_{k+1}(1)=1$, it follows that
$\delta([1])=0$. Let 
\begin{displaymath}
Y:= (1\otimes p_{n-2}\otimes 1_{k})(1\otimes 1_{n-2} \otimes p_{k-1}\otimes
1)\tilde{X}+ 1-1\otimes p_{n-2}\otimes p_{k-1} \otimes 1  .
\end{displaymath}
Since $1\otimes p_{n-2}\otimes 1  =
1_{\{1\}}(\chi_{\omega_{k}}(u_{n-1,1}^{*}u_{n-1,1})$ and $
 1\otimes 1_{n-2} \otimes p_{k-1}  = 1_{\{1\}}(\tilde{X}^{*}\tilde{X})$ it
follows that the operator 
 $Y \in C(S_{q}^{n,2,k+1})$. Also 
\begin{displaymath}
 Y= t_{1}\otimes p_{n-2}\otimes p_{k-1}\otimes S^{*} + 1-1\otimes p_{n-2}\otimes
p_{k-1}\otimes 1.
\end{displaymath}
Note that $Y$ is an isometry such that $\sigma_{k+1}(Y)=u_{k}$. One has
$\sigma_{k+1}(v_{k+1})=v_{k}$. Note that $Y$ and $v_{k+1}$ commute. Hence by
corollary \ref{index}, it follows that the image of $\delta$ is the subgroup
generated by $[v_{k+1}(1-YY^{*})+YY^{*}]=[V_{1}\otimes p_{k}+1-1\otimes
p_{k}]$. 

Thus the above computation with the six term sequence implies that
$K_{0}(C(S_{q}^{n,2,k+1}))$ is isomorphic to $\mathbb{Z}^{2}$ and is generated
by $P(u_{1},v_{1})\otimes p_{k}=P(u_{k},v_{k})$ and $[1]$ and
$K_{1}(C(S_{q}^{n,2,k+1}))$ is isomorphic to $\mathbb{Z}^{2}$ and is generated
by $[V_{k+1}]$ and $[U_{1}\otimes p_{k}+1-1\otimes p_{k}]=[U_{k+1}]$. This
completes the proof. \hfill $\Box$

\section{$K$ groups of $C(S_{q}^{n,2})$}
In this section we compute the $K$ groups of $C(S_{q}^{n,2})$. We start with a
few  observations.
\begin{lmma}
In the permutation group $S_{n}$ one has
$\omega_{n-2,1}\omega_{n-1,1}=\omega_{n-1,1}\omega_{n-1,2}$. 
\end{lmma}
\textit{Proof.} First note that $s_{i}s_{i+1}s_{i}= s_{i+1}s_{i}s_{i+1}$ and
$s_is_{j}=s_{j}s_{i}$ if $|i-j| \geq 2$. Hence one has
$\omega_{n-1,k}\omega_{n-1,1}=\omega_{n-1,k+1}\omega_{n-1,1}s_{k+1}$. Now the
result follows by induction on $k$. \hfill $\Box$

 We denote the representation $\chi_{\omega_{n-1,1}}*\pi_{\omega_{n-1,2}}$ by
$\tilde{\chi}_{\omega_{n}}$. Since $\omega_{n-1,1}\omega_{n-1,2}$ is a reduced
expression for $\omega_{n}$ it follows that the representations
$\tilde{\chi}_{\omega_{n}}$ and $\chi_{\omega_{n}}$ are equivalent. Let $U$ be a
unitary such that $U\chi_{\omega_{n}}(.)U^{*}=\tilde{\chi}_{\omega_{n}}(.)$. It
is clear that $\tilde{\chi}_{\omega_{n}}(C(S_{q}^{n,2})) \subset
C(\mathbb{T}^{m})\otimes \tau \otimes \tau^{\otimes \ell(\omega_{n-1})}$. Let
$\tilde{\sigma}_{n}$ denote the restriction of $1\otimes \sigma \otimes
1^{\otimes(2(n-2)}$ to $\tilde{\chi}_{\omega_{n}}(C(S_{q}^{n,2}))$. Since
$\tilde{\sigma}_{n}(\tilde{\chi}_{\omega_{n}}(u_{ij}))=
\chi_{\omega_{n-1}}(u_{ij})$ one has the following commutative diagram

\begin{tikzpicture}\label{com}[description/.style={fill=white,inner sep=2pt}]
    \matrix (m) [matrix of math nodes, row sep=6em,
    column sep=2.5em, text height=1.5ex, text depth=0.25ex]
    { & &\chi_{\omega_{n}}(C(S_{q}^{n,2})) & &
\tilde{\chi}_{\omega_{n}}(C(S_{q}^{n,2})) \\
      & & & C(S_{q}^{n,2,n-1}) & \\ };
      \path[->,font=\scriptsize]
    (m-1-3) edge node[auto] {$ U(.)U^{*} $} (m-1-5)
            edge node[below] {$ \sigma_{n} $} (m-2-4)
    (m-1-5) edge node[auto] {$ \widetilde{\sigma}_{n}$} (m-2-4);
\end{tikzpicture}
\begin{lmma}
\label{coisometry}
 There exists a coisometry $X \in \chi_{\omega_{n}}(C(S_{q}^{n,2}))$ such that
$\sigma_{n}(X)=V_{n-1}$ and
$X^{*}X=1-1_{\{1\}}(\chi_{\omega_{n}}(u_{n1}^{*}u_{n1}))$. 
\end{lmma} 
\textit{Proof.} By the above commutative diagram,  it is enough to show that
there exists a coisometry $\widetilde{X} \in
\widetilde{\chi}_{\omega_{n}}(C(S_{q}^{n,2}))$ such that
$\widetilde{\sigma_{n}}(X)=V_{n-1}$ and
$X^{*}X=1-1_{\{1\}}(\widetilde{\chi}_{\omega_{n}}(u_{n1}^{*}u_{n1})$.
Now note that
$\widetilde{\chi}_{\omega_{n}}(u_{n-1,1}^{*}u_{n-1,1})-q^{2}u_{n1}u_{n1})=
1\otimes 1 \otimes \underbrace{q^{2N}\otimes q^{2N}\otimes \cdots q^{2N}}_{(n-2)
times} \otimes 1_{n-2}$. Hence the projection $1\otimes 1 \otimes p_{n-2}
\otimes
1_{n-2}=1_{\{1\}}(\widetilde{\chi}_{\omega_{n}}(u_{n-1,1}^{*}u_{n-1,1}-q^{2}u_{
n1}^{*}u_{n1}))$ is in the algebra
$\widetilde{\chi}_{\omega_{n}}(C(S_{q}^{n,2}))$. Now let $Y:=(1\otimes 1 \otimes
p_{n-2} \otimes 1_{n-2})\widetilde{\chi}_{\omega_{n}}(u_{n-1,1})$. Then
$Y:=t_{2} \otimes \sqrt{1-q^{2N+2}}S \otimes p_{n-2} \otimes 1_{n-2}$. Hence the
operator $Z:=t_{2} \otimes S \otimes p_{n-2} \otimes 1_{n-2}$ is in the algebra
$\widetilde{\chi}_{\omega_{n}}(C(S_{q}^{n,2}))$. Now let
$\widetilde{X}:=Z+1-ZZ^{*}$. Then $\widetilde{X}$ is a coisometry such  that
$\widetilde{\sigma}_{n}(\widetilde{X})=V_{n-1}$ and
$\widetilde{X}^{*}\widetilde{X}=1-1\otimes p_{n-1} \otimes 1_{n-2}$ which is
$1-1_{\{1\}}(\widetilde{\chi}_{\omega_{n}}(u_{n1}^{*}u_{n1}))$. This completes
the proof. \hfill $\Box$

Observe that the operator $\tilde{Z}:=t_{1}\otimes 1_{n-1} \otimes
\underbrace{q^{N}\otimes q^{N} \otimes \cdots q^{N}}_{(n-2)times} \otimes S^{*}$
 lies in the algebra $C(S_{q}^{n,2,n})$ since the difference
$\tilde{Z}-\chi_{\omega_{n}}(u_{n,2})$ lies in the ideal
$C(S_{q}^{n,2,1})\otimes \mathcal{K}^{\otimes(n-1)})$. Let
 $Z:=1_{\{1\}}(\tilde{Z}^{*}\tilde{Z})\tilde{Z}$ and $Y_{n}:=Z+1-Z^{*}Z$. Then 
\begin{eqnarray}
\label{defn of Z_{n}}
Z_{n}&=&t_{1} \otimes 1_{n-2} \otimes p_{n-2} \otimes S^{*}, \\
\label{defn of Y_{n}}
Y_{n}&=&t_{1} \otimes 1_{n-2} \otimes p_{n-2} \otimes S^{*} + 1-1\otimes 1_{n-2}
\otimes p_{n-2} \otimes 1.
\end{eqnarray}
 Hence $Y$ is an isometry and
$YY^{*}=1-1_{\{1\}}(\chi_{\omega_{n}}(u_{n1}^{*}u_{n1}))$. Let $X$ be a
coisometry in $C(S_{q}^{n,2,n})$ such that $\sigma_{n}(X)=v_{n-1}$ and
$X^{*}X:=1-1_{\{1\}}(\chi_{\omega_{n}}(u_{n1}^{*}u_{n1}))$. The existence of
such an $X$ was shown in lemma \ref{coisometry}. Then $XY$ is a unitary. 
\begin{ppsn}
\label{K group computation}
 The $K$ groups $K_{0}(C(S_{q}^{n,2})$ and $K_{1}(C(S_{q}^{n,2})$ are both
isomorphic to $\mathbb{Z}^{2}$. In particular we have the following.
\begin{enumerate}
 \item The projections $[1]$ and $P(u_{n},v_{n})$ generate
$K_{0}(C(S_{q}^{n,2}))$.
 \item The unitaries $U_{n}$ and $XY_{n}$ generate $K_{1}(C(S_{q}^{n,2}))$ where
$X$ is a coisometry in $C(S_{q}^{n,2})$ such that $\sigma_{n}(X)=V_{n-1}$ and
$X^{*}X=1-1_{\{1\}}(u_{n1}^{*}u_{n1})$ and $Y_{n}$ is as in equation \ref{defn
of Y_{n}}
\end{enumerate}
\end{ppsn}
\textit{Proof.} By proposition \ref{sheu}, we have the following exact sequence.
\begin{displaymath}
0 \longrightarrow C(S_{q}^{n,2,1})\otimes \mathcal{K}^{\otimes(n-1)}
\longrightarrow C(S_{q}^{n,2,n}) \stackrel{\sigma_{n}} \longrightarrow
C(S_{q}^{n,2,n-1}) \longrightarrow 0.
\end{displaymath}
which gives rise to the following six term sequence in $K$ theory.
\begin{equation*}
\def\labelstyle{\scriptstyle}
\xymatrix@C=25pt@R=20pt{
K_0(C(S_{q}^{n,2,1})\otimes \mathcal{K}^{\otimes n-1}) \ar[r]&
K_0(C(S_{q}^{n,2,n}))\ar[r]^{K_0(\sigma_{n})}& K_0(C(S_{q}^{n,2,k}))
\ar[d]_{\delta} \\
K_1(C(S_{q}^{n,2,n-1})) \ar[u]^{\partial} &
K_1(C(S_{q}^{n,2,n}))\ar[l]^{K_{1}(\sigma_{n})}& K_1(C(S_{q}^{n,2,1})\otimes
\mathcal{K}^{\otimes n-1})\ar[l] 
}
\end{equation*}
Now we compute $\partial$ and $\delta$ to compute the six term sequence. First
note that since $[U_{n-1}]$ and $[V_{n-1}]$ generate
$K_{1}(C(S_{q}^{n,2,n-1}))$, it follows that $[U_{n-1}]$ and $[V_{n-1}U_{n-1}]$
generate $K_{1}(C(S_{q}^{n,2,n-1}))$. As $XY_{n}$ is a unitary for which
$\sigma_{n}(XY_{n})=V_{n-1}U_{n-1}$, it follows that
$\partial([V_{n-1}U_{n-1}])=0$. Next $Y_{n}$ is an isometry for which
$\sigma_{n}(Y_{n})=U_{n-1}$. Hence $\partial([U_{n-1}])=[1-Y^{*}Y]-[1-YY^{*}]$.
Thus $\partial([U_{n-1}])=-[1\otimes 1_{n-2}\otimes p_{n-1}]$.

Now we compute $\delta$. Since $\sigma_{n}(1)=1$, it follows that
$\delta([1])=0$. Now one observes that $p_{n-2} \otimes
S^{*}\pi_{\omega_{n-1,1}}(u_{j1})=0$ if $j>1$. Hence
$Z_{n}\chi_{\omega_{n}}(u_{n-1,1})=t_{1}t_{2} \otimes p_{n-2} \otimes p_{n-2}
\otimes \sqrt{1-q^{2N+2}}$ where $Z_{n}$ is as defined in \ref{defn of Z_{n}}.
Thus the operator $R_{n}:=t_{1}t_{2}\otimes p_{n-2} \otimes p_{n-2} \otimes 1$
lies in the algebra $C(S_{q}^{n,2,n})$ as the difference
$R_{n}-Z_{n}\chi_{\omega_{n}}(u_{n-1,1})$ lies in the ideal
$C(\mathbb{T}^{2})\otimes \mathcal{K}^{\otimes(2n-3))}$. Hence projection
$1\otimes p_{n-2} \otimes p_{n-2} \otimes 1$ lies in the algebra
$C(S_{q}^{n,2,n})$. Now define
\begin{displaymath}
\begin{array}{ll}
 S_{n}:=R_{n}+1-R_{n}R_{n}^{*},\\
 T_{n}:=(1\otimes p_{n-2} \otimes p_{n-2} \otimes 1)Z_{n}+ 1-1\otimes
p_{n-2}\otimes p_{n-2} \otimes 1.
\end{array}
\end{displaymath}
Then $S_{n}$ is a unitary and $T_{n}$ is an isometry such that
$\sigma_{n}(S_{n})=u_{n-1}v_{n-1}$ and $\sigma_{n}(T_{n})=u_{n-1}$. Moreover
$S_{n}$ and $T_{n}$ commute. Now note that
$P(u_{n-1},v_{n-1})=P(u_{n-1},u_{n-1}v_{n-1})$. Hence by corollary \ref{index},
it follows that the image of $\delta$ is the subgroup generated by
$S_{n}(1-T_{n}T_{n}^{*})+T_{n}T_{n}^{*}$ in $K_{1}(C(S_{q}^{n,2,1})\otimes
\mathcal{K}^{\otimes(n-1)})$. Now 
\begin{displaymath}
 S_{n}(1-T_{n}T_{n}^{*})+T_{n}T_{n}^{*}=t_{1}t_{2}\otimes p_{n-2}\otimes
p_{n-1}+1-1\otimes p_{n-2} \otimes p_{n-1}.
\end{displaymath}
Since $1\otimes p_{n-2}$ is a trivial  in $K_{0}(C(S_{q}^{2n-3}))$ it follows
that the unitary $t_{1}\otimes p_{n-2} + 1-1\otimes p_{n-2}$ is trivial in
$K_{1}(C(S_{q}^{n,2,1}))=K_{1}(C(\mathbb{T})\otimes C(S_{q}^{2n-3})$. Hence one
has  $[S_{n}(1-T_{n}T_{n}^{*})+T_{n}T_{n}^{*}]=[V_{1}\otimes p_{n-1}+1-1\otimes
p_{n-1}]$ in $K_{1}(C(S_{q}^{n,2,1})\otimes \mathcal{K}^{\otimes(n-1)})$.

Thus the above computation with the exactness of the six term sequence completes
the proof. \hfill $\Box$
\section{K groups of quantum SU(3)}
In this section we show that for $n=3$ the unitary $XY_{n}$ in proposition
\ref{K group computation} can be replaced by the fundamental $3\times 3$ matrix
$(u_{ij})$ of $C(SU_{q}(3))$. First note that for $n=3$ we have
$C(S_{q}^{n,2})=C(SU_{q}(3))$ since $C(SU_{q}(1))=\mathbb{C}$. The embedding $SU_q(1) \subseteq SU_q(3)$ is given by the counit. Hence the quotient
$C(SU_q(3)/SU_q(1))$ becomes isomorphic with $C(SU_q(3))$. The algebra
$C(S_{q}^{3,2,1})$ is denoted $C(U_{q}(2))$ in \cite{Sh1}. Then
$C(U_{q}(2))=C(\mathbb{T})\otimes C(SU_{q}(2))$. Let $ev_{1}:C(\mathbb{T})\to
\mathbb{C}$ be the evaluation at the point '1'. Then $\phi=(ev_{1}\otimes
1)\sigma_{2}\sigma_{3}$ where $\phi:C(SU_{q}(3))\to C(SU_{q}(2))$ is the
subgroup homomorphism defined in equation \ref{subgroup}.
\begin{ppsn}
The $K$ group $K_{1}(C(SU_{q}(3))$ is isomorphic to $\mathbb{Z}^{2}$ generated
by the unitary $U_{3}:=t_{1}\otimes p\otimes p+1-1\otimes p\otimes p$ and the
fundamental unitary $U=(u_{ij})$ 
\end{ppsn}
\textit{Proof.} By proposition \ref{K group computation}, we know that
$K_{1}(C(SU_{q}(3))$ is isomorphic to $\mathbb{Z}^{2}$ and is generated by
$[U_{3}]$ and $[XY_{3}]$ where $X$ is a coisometry such that
$\sigma_{3}(X)=V_{2}$ and
$X^{*}X=1-1_{\{1\}}(\chi_{\omega_{3}}(u_{31}^{*}u_{31}))$. Now observe that
$\phi(X)=t_{2}\otimes p+1-1\otimes p$ and $\phi(Y_{3})=1$. Hence
$\phi(XY_{3})=t_{2}\otimes p+1-1\otimes p$. Also note that $\phi(U_{3})=0$ and
$\phi(U)=\begin{bmatrix}
                                             u & 0 \\                           
                            0 & 1

\end{bmatrix}$ where $u$ denote the fundamental unitary of $C(SU_{q}(2))$. Since
$K_{1}(C(SU_{q}(2))$ is isomorphic to $\mathbb{Z}$ the proof is complete if we 
show that $t_{2}\otimes p+1-1\otimes p$ and $[u]$ represents the same element in
$K_{1}(C(SU_{q}(2))$ which we do in the next lemma.\hfill$\Box$

We denote the $2\times 2$ fundamental unitary $u=(u_{ij})$ of $C(SU_{q}(2))$ by
$u_{q}$. Consider the representation $\chi_{s_{1}}:C(SU_{q}(2))\to
B(\ell^{2}(\mathbb{Z})\otimes \ell^{2}(\mathbb{N}))$. We let the unitary $t$ act
on $\ell^{2}(\mathbb{Z})$ as the right shift i.e $te_{n}=e_{n+1}$. Let
$\{e_{n,m}: n \in \mathbb{Z}, m \in \mathbb{N}\}$ be the standard orthonormal
basis for the Hilbert space $\ell^{2}(\mathbb{Z}) \otimes \ell^{2}(\mathbb{N})$.
For an integer $k$, denote the orthogonal projection onto the closed subspace
spanned by $\{e_{n,m}:n+m \leq k\}$ by $P_{k}$ and set $F_{k}:=2P_{k}-1$. Note
that $F_{k}$ is a selfadjoint unitary.
\begin{ppsn}
For any integer $k$, the triple $(\chi_{s_{1}}, \ell^{2}(\mathbb{Z}) \otimes
\ell^{2}(\mathbb{N}), F_{k})$ is an odd Fredholm module for $C(SU_{q}(2))$ and
we have the pairing
\begin{enumerate}
\item $ \langle [u_{q}],F_{k} \rangle =-1$
\item $ \langle t \otimes p + 1-1 \otimes p, F_{k} \rangle =-1$ where
$p=1-S^{*}S$.
\end{enumerate}
\end{ppsn}
\textit{Proof.} It is not difficult to show that $C(SU_{q}(2))$ is generated by
$t \otimes S$ and $t \otimes p$. Now it is easy to see that $[t \otimes
S,P_{k}]=0$ and $[t \otimes p,P_{k}]$ is a finite rank operator. Hence the
triple $(\chi_{s_{1}}, \ell^{2}(\mathbb{Z}) \otimes \ell^{2}(\mathbb{N}),
F_{k})$ is an odd Fredholm module for $C(SU_{q}(2))$. Since $C(SU_{q}(2))$ is
generated by $t \otimes S$ and $t \otimes p$ it follows that $u_{p} \in
C(SU_{q}(2))$ for every $p > 0$. Also as $p \to 0$, $u_{p}$ approaches to $u$ in
norm where $u$ is given by
\begin{displaymath}
\begin{array}{ll}
u:= & \left ( \begin{array}{ll}
                  t \otimes S & ~~~0\\
                 \bar{t} \otimes p & \bar{t} \otimes S^{*}  
                  \end{array} \right)
\end{array}
\end{displaymath}
Hence $[u_{q}]=[u]$ in $K_{1}(C(SU_{q}(2)))$. It is easy to check  the following
\begin{align*}
\langle [u],F_{k} \rangle&=-1 \\
\langle [t \otimes p+ 1-1\otimes p],F_{k} \rangle&=-1.
\end{align*}
This completes the proof. \hfill $\Box$

\bibliography{references}

\def\cprime{$'$} \def\cprime{$'$}
\begin{thebibliography}{10}

\bibitem{Bla}
B.~Blackadar.
\newblock K-theory for operator algebras.
\newblock Springer Verlag,Newyork, 1987.

\bibitem{CPR_arxiv}
Alan~L. Carey, John Phillips, and Adam Rennie.
\newblock Twisted cyclic theory and an index theory for the gauge invariant kms
  state on cuntz algebras.
\newblock arXiv:0801.4605.

\bibitem{PsPal_SU_q_2}
Partha~Sarathi Chakraborty and Arupkumar Pal.
\newblock Equivariant spectral triples on the quantum {${\rm SU}(2)$} group.
\newblock {\em $K$-Theory}, 28(2):107--126, 2003.

\bibitem{CP_Tran}
Partha~Sarathi Chakraborty and Arupkumar Pal.
\newblock Equivariant spectral triples and poincaré duality for $su_q(2)$.
\newblock {\em Trans. Amer. Math. Soc.}, 2010.

\bibitem{Con}
A.~Connes.
\newblock An analogue of the thom isomorphism for a crossed product of a
  ${C}^{*}$ algebra by an action of $\mathbb{R}$.
\newblock {\em Advances in Mathematics}, 39(1):31--55, 1981.

\bibitem{ConIHES}
Alain Connes.
\newblock Noncommutative differential geometry.
\newblock {\em Inst. Hautes \'Etudes Sci. Publ. Math.}, (62):257--360, 1985.

\bibitem{ConSU_q2}
Alain Connes.
\newblock Cyclic cohomology, quantum group symmetries and the local index
  formula for {${\rm SU}\sb q(2)$}.
\newblock {\em J. Inst. Math. Jussieu}, 3(1):17--68, 2004.

\bibitem{CM_typeIII}
Alain Connes and Henri Moscovici.
\newblock Type {III} and spectral triples.
\newblock In {\em Traces in number theory, geometry and quantum fields},
  Aspects Math., E38, pages 57--71. Friedr. Vieweg, Wiesbaden, 2008.

\bibitem{Le-So}
Serge Levendorski{\u\i} and Yan Soibelman.
\newblock Algebras of functions on compact quantum groups, {S}chubert cells and
  quantum tori.
\newblock {\em Comm. Math. Phys.}, 139(1):141--170, 1991.

\bibitem{McClanahan}
Kevin McClanahan.
\newblock {$C^*$}-algebras generated by elements of a unitary matrix.
\newblock {\em J. Funct. Anal.}, 107(2):439--457, 1992.

\bibitem{Nagy}
Gabriel Nagy.
\newblock Bivariant {$K$}-theories for {$C^*$}-algebras.
\newblock {\em $K$-Theory}, 19(1):47--108, 2000.

\bibitem{NT_arxiv}
Sergey Neshveyev and Lars Tuset.
\newblock The dirac operator on compact quantum groups.
\newblock arXiv:math/0703161.

\bibitem{VP}
G.B. Podkolzin and L.I. Vainerman.
\newblock Quantum stiefel manifold and double cosets of quantum unitary group.
\newblock {\em Pacific Journal of Mathematics}, 188(1), 1999.

\bibitem{Pod}
Piotr Podle{\'s}.
\newblock Symmetries of quantum spaces. {S}ubgroups and quotient spaces of
  quantum {${\rm SU}(2)$} and {${\rm SO}(3)$} groups.
\newblock {\em Comm. Math. Phys.}, 170(1):1--20, 1995.

\bibitem{Sh-pacific}
Albert Jeu-Liang Sheu.
\newblock The structure of twisted {${\rm SU}(3)$} groups.
\newblock {\em Pacific J. Math.}, 151(2):307--315, 1991.

\bibitem{Sh1}
Albert~J.L. Sheu.
\newblock Compact quantum groups and groupoid ${C}^\ast$ algebras.
\newblock {\em Journal of Functional Analysis}, 144:371--393, 1997.

\bibitem{So1}
Ya.~S. So{\u\i}bel{\cprime}man.
\newblock Algebra of functions on a compact quantum group and its
  representations.
\newblock {\em Algebra i Analiz}, 2(1):190--212, 1990.

\bibitem{Vak-Soi}
L.~L. Vaksman and Ya.~S. So{\u\i}bel{\cprime}man.
\newblock An algebra of functions on the quantum group {${\rm SU}(2)$}.
\newblock {\em Funktsional. Anal. i Prilozhen.}, 22(3):1--14, 96, 1988.

\bibitem{Wor_twisted}
S.~L. Woronowicz.
\newblock Twisted {${\rm SU}(2)$} group. {A}n example of a noncommutative
  differential calculus.
\newblock {\em Publ. Res. Inst. Math. Sci.}, 23(1):117--181, 1987.

\bibitem{Wor_Tannaka}
S.~L. Woronowicz.
\newblock Tannaka-{K}re\u\i n duality for compact matrix pseudogroups.
  {T}wisted {${\rm SU}(N)$} groups.
\newblock {\em Invent. Math.}, 93(1):35--76, 1988.

\end{thebibliography}
\bibliographystyle{plain}
\end{document}